\documentclass[12pt]{article}
\usepackage{latexsym,color,amsmath,amsthm,amssymb,amscd,amsfonts}

\newtheorem{Theorem 3.9}{Theorem 3.9}
\newtheorem{Corollary}{Corollary}[section]
\newtheorem{Non Symmetric Example}{Non Symmetric Example}[section]

\newtheorem{Fact}{Fact}

\newtheorem{lemma}{Lemma}[section]

\newtheorem{theorem}[lemma]{Theorem}

\newtheorem{corollary}[lemma]{Corollary}
\newtheorem{thm}{Theorem}

\newtheorem{rem}[thm]{Remark}

\newtheorem*{remark*}{Remark}

\def\L{{\cal L}}

\def\R{{\mathbb R}}
\def\N{{\mathbb N}}
\def\grad{\nabla}

\def\qed{\hfill $\vcenter{\hrule height .3mm
\hbox {\vrule width .3mm height 2.1mm \kern 2mm \vrule width .3mm
height 2.1mm} \hrule height .3mm}$ \bigskip}

\begin{document}

{\setlength{\baselineskip}%
{1.5\baselineskip}

\makebox[\textwidth]{%
 \Huge  Santal\'{o} region of a log-concave function }\par

\par}

\begin{abstract}
In this paper  we defined the {\it Santal\'{o}-region} and the {\it Floating body} of a log-concave function on $\R^n$. We
then study their properties. Our main result is that any relation of {\it Floating body} and {\it Santal\'o region} of a
convex body  is translated to a relation of {\it Floating body} and {\it Santal\'o region} of an even log-concave function.

\end{abstract}

\footnote{T.Weissblat   Ph.D student Tel Aviv University 69978 Israel.}

\newpage

\section{Introduction}

 In $\R^n$ we fix a
scalar product $\langle \cdot , \cdot \rangle$. A convex body in
$\R^n$ is a compact convex set which includes $0$ in its interior.
We denote by $B_{2}^{n}$ the {\it Euclidean} unit ball of radius $1$ in
$\R^{n}$, namely $B_2^n = \{ x\in \R^n : \sum_i x_i^2 \le 1\}.$ We
let the boundary of the ball, namely the unit sphere, be denoted by
$S^{n-1}$. We call a body centrally symmetric if $K= -K$.

Let $K$ be a convex body in $\R^{n}$. For $x\in int(K)$, the interior of $K$, let $K^{x}$ be the polar body of
$K$ with respect to $x$, i.e \[ K^{x}= \{x+y:y \in \R^{n},\langle y,z-x \rangle\leq1 \ \ \forall z\in
K\}=x+(K-x)^{\circ}\]  where $ K^{\circ} := \{y \in
\R^{n} : \langle y,z \rangle\leq1 \ \ \forall z\in K \}$, is the standard polar of a body.

The {\it Blaschke-Santal\'o inequality}  asserts
that the quantity  $\inf_x Vol(K)Vol(K^{x})$ is maximal for balls
and ellipsoids. Moreover, it is well known that there exists a
unique $x_{0}\in int(K)$ on which this infimum is attained.  This unique $x_0$ is called the
{\it Santal\'{o} point} of $K$.

The above gives rise to a natural geometric construction which was
invented by E.~Lutwak and consists of all the points which are
``witnesses'' of the validity of the {\it Blaschke-Santal\'o inequality}.
This is called the {\it Santal\'o-region} of $K$ (with constant 1) and is
defined by
$ S(K,1)=\{x\in K: Vol(K)Vol(K^{x})\leq Vol(B_{2}^{n})^{2} \}$. More generally, one defines the {\it Santal\'o-region} of $K$ with
constant $t$, denoted $S(K, t)$, by
\begin{eqnarray}\label{santalo region of K}
S(K,t)=\{x\in K: Vol(K)Vol(K^{x})\leq t Vol(B_{2}^{n})^{2} \}.
\end{eqnarray}
It follows from the {\it Blaschke-Santal\'o inequality}  that the {\it Santal\'o point} $x_{0}\in S(K,1)$, and that
$S(K,1)=\{x_{0}\} $ if and only if $K $ is an ellipsoid.

We turn to the definition of the {\it Floating body}. For
$0<\lambda<\frac{1}{2}$ and $\theta\in S^{n-1}$ we define $H_{\lambda,\theta}:=\{x\in \R^{n}:\langle x,\theta \rangle = a(\theta) \}$, where
$a(\theta)\in\R$ is such that
$Vol(H_{\lambda,\theta}^{+}\cap K)=\lambda Vol(K)$,
and $H_{\lambda,\theta}^{+}=\{x\in \R^{n}: \langle x,\theta \rangle \geq a(\theta) \}$. Meaning,
$H_{\lambda,\theta}$ is a hyperplane in the direction $\theta$ that 'cuts' $\lambda Vol(K)$ volume from $K$. In what follows we will use the notation $H_{\lambda,\theta}^{-}$, which means
$H_{\lambda,\theta}^{-}=\{x\in \R^{n}: \langle x,\theta
\rangle \leq a(\theta) \}$.

Let $K$ be a convex body in $\R^{n}$ and $0<\lambda<\frac{1}{2}$. We define

\begin{eqnarray}\label{floating body of K}
F(K,\lambda)=\bigcap_{\theta\in S^{n-1}}H_{\lambda,\theta}^{-}.
\end{eqnarray}

This body is called the {\it Floating body} of $K$, and was first
considered by Dupin in 1822. A vast literature exists on the topic
of  Floating bodies, see e.g. \cite{mr},\cite{mr2},\cite{mw}.
 We list some known results about  Floating bodies
which exemplify their importance in Section \ref{backg}. For us here it is important to note that
 $F(K,\lambda)$ is a closed convex set which is contained in
$K$ (possibly empty). In the centrally symmetric case it is known that the family supporting hyperplanes
touching $F(K, \lambda)$ at its boundary point are exactly all of $H_{\lambda, \theta}$ above, so that $F(K,
\lambda)$ is the envelope of these hyperplanes. We state this result precisely in Section \ref{backg} below. In
fact, we state there a slightly more general theorem valid for measures, not just bodies.

From convex bodies we turn to log-concave functions. A function $f:\R^{n}\rightarrow [0,\infty)$ is called
log-concave
 if $$\lambda log f(x)+(1-\lambda)log f(y)\leq log f(\lambda
x + (1-\lambda)y)$$ for any $0<\lambda<1 $ and $x,y\in supp(f),$ where $supp(f)$ denotes the support of the
function $f$, $ supp (f) = \{x\in\R^{n}: f(x)>0 \}$. We say that a function is even if $f(x) = f(-x)$ for all
$x\in\R^{n}$.
 The idea of transferring results about convex bodies to results about log-concave functions is not new, and
its roots can be found in \cite{ba86}. We will elaborate on the importance of this procedure, and give examples
for it, in Section \ref{backg} below. Let us just mention that although, so far, it seems that almost every
geometric result on convex bodies has a counterpart which is valid for log-concave functions, and although
several methods for obtaining such results have been found, there is so far no general procedure for doing this,
and this is more of a meta-mathematical fact that an actual mathematical truth.

Motivated by \cite{AKM},  we define the dual of a log-concave function $f$ by
\begin{eqnarray*}
 f^{\circ}(x)=inf_{y\in \R^{n}}[\frac{e^{-\langle x,y
\rangle}}{f(y)}].
\end{eqnarray*}
 We explain the motivation for this definition and equivalent formulations in Section \ref{backg}. It turns out
that a {\it Santal\'o type inequality} holds also for functions. Again, one should choose correctly the center of the
function. Thus, for a log-concave function $f$, we let $f_a$ stand for the translation of $f$ by the vector $a$,
namely $f_{a}(x)= f(x-a)$.

When ``volume'' is exchanged with ``integral'' and ``an ellipsoid''
is exchanged with ``a Gaussian function'', we have a {\it Santal\'o type
inequality} for functions.

\begin{theorem}\cite{AKM}\label{vitali} Let $f:\R^{n}\rightarrow[0,\infty)$ be
an integrable function such that $0<\int f<\infty$. Then, for some
vector $x_{0}$, one has
$$ \int_{\R^{n}}  f_{x_{0}}\int_{\R^{n}} (f_{x_{0}})^{\circ}\leq
(2\pi)^{n}.$$ If $f$ is log-concave, one may choose
$x_{0}=\frac{\int xf(x)}{\int f(x)}$, the center mass of $f$.
The minimum over $x_{0}$ of the left side product equals
$(2\pi)^{n}$ if and only if $f$ is proportional to Gaussian, namely
$f=ce^{-\langle Ax,x\rangle}$, where $0<c$ and $A$ is a
positive-definite matrix.
\end{theorem}

 As in the case of bodies, we call the point $x_{0}$ which minimizes $a\rightarrow \int (f_{a})^{\circ}dx$
the {\it Santal\'o point} of $f$. The case of even functions $f$ (in which the {\it Santal\'o point} is $0$) was proven much
earlier in  \cite{ballphd}. It is natural to define for a log-concave function
its {\it Santal\'o region} (with parameter $t$) by
\begin{eqnarray}\label{santalo region of f}
S(f,t)=\{a\in \R^{n}: \int f_{a} \int (f_{a})^{\circ}\leq (2\pi)^{n}t
\}.
\end{eqnarray}

\noindent  It follows from Theorem \ref{vitali} that the {\it Santal\'o point}
$x_{0} \in S(f,1)$, and that $S(f,1)=\{x_{0}\}$ if and only
if $f$ is proportional to Gaussian.

Similarly to the {\it Floating body} of a convex body, one may define the
{\it Floating body} of a log concave function $f$. Let $0<\lambda<\frac{1}{2}$, and   consider $H_{\lambda,\theta}=\{x\in \R^{n}:\langle
x,\theta \rangle = a(\theta) \}$, where $a(\theta)$ satisfies
$\int_{H_{\lambda,\theta}^{+}}f dx= \lambda \int fdx$. We define the {\it Floating body} of $f$
\begin{eqnarray}
F(f,\lambda)= \bigcap_{\theta \in S^{n-1}}H_{\lambda,\theta}^{-}.
\end{eqnarray}

Notice that if $f$ is an even function, we have that its {\it Floating body} is non-empty set for any
$0<\lambda<\frac{1}{2}$, since $0\in F(f,\lambda)$. Moreover, it was proved by Meyer-Reisner that for every even log-concave function $f:\R^{n}\rightarrow[0,\infty)$,  the family
supporting hyperplanes touching $F(f, \lambda)$ at its boundary point are exactly all of $H_{\lambda, \theta}$
above, so that $F(f, \lambda)$ is the envelope of these hyperplanes (section \ref{backg}).

The main result of this paper is that any 'relation' of {\it Floating body} and {\it Santal\'o region} of a
convex body is translated to a relation of {\it Floating body} and {\it Santal\'o region} of an even log-concave function, namely

\begin{thm} \label{Thm}
Let $0<\lambda<\frac{1}{2}$ and $0<d$ such that
 $F(K,\lambda) \subseteq S(K,d)$ for every  centrally
symmetric convex body K. Then $$F(f,\lambda) \subseteq S(f,d) \textrm{  for every  even log-concave
function } f:\R^{n} \rightarrow [0,\infty).$$
\end{thm}

 In the case of convex bodies, the following
inclusion  holds

\begin{theorem}\cite{mw}\label{Meyer}
Let $K$ be a convex body in $\R^{n}$. Then

  $$F(K,\lambda) \subseteq
S(K,\frac{1}{4\lambda(1-\lambda)}) \textrm{ for any }0<\lambda <1/2.$$
\end{theorem}

From Theorem \ref{Thm} and Theorem \ref{Meyer} we get
\begin{corollary}
For every even log-concave function $f:\R^{n}\rightarrow [0,\infty)$
and any $0<\lambda<\frac{1}{2}$ we have that
$$F(f,\lambda)\subseteq S(f,\frac{1}{4\lambda(1-\lambda)}).$$
\end{corollary}

It is not known whether - in some form - the {\it Floating body} and the
{\it Santal\'o region} of a convex body are isomorphic. However, we can
show that any inclusion in the direction opposite to that of Theorem
\ref{Meyer} will transfer immediately to the realm of functions.
That is, the following ``converse'' of Theorem \ref{Thm} holds.

\begin{thm}\label{reverse} Let $0<\lambda<\frac{1}{2}$ and
$0<d$ such that $S(K,d) \subseteq F(K,\lambda)$ for every  centrally
symmetric convex body K. Then
$$S(f,d) \subseteq F(f,\lambda) \textrm{ for every even log-concave
function}  f:\R^{n} \rightarrow [0,\infty).$$
\end{thm}

The rest of the paper is organized as follows. In section \ref{backg} we explore in depth the notions defined
above. We quote the relevant known theorems, and discuss several other notions which we will need in the
proofs of our main theorems, such as $s$-concave functions. In section \ref{TSR} we explore the properties of
the {\it Santal\'o region} of a function. In Section \ref{bigsection} we prove Theorem \ref{Thm} and Theorem
\ref{reverse}.

\newpage

\section{Background}\label{backg}

The following theorem, namely Theorem \ref{IN}, is the main ingredient
for several important things. First,   the {\it Blaschke-Santal\'o inequality }    is an immediate result of it.
Second, it gives us a connection (Theorem \ref{Meyer}) between the {\it Floating body} and the  {\it Santal\'o region} which  related to a convex body, which afterwards we translate (Theorem \ref{Thm})  to a connection between the
{\it Floating body} and the {\it Santal\'o region}
 related to an even log-concave function.

\begin{theorem}\cite{mp1}\label{IN}
Let $K$ be a convex body in $\R^{n}$ and let $H=\{x\in\R^{n}:\langle
x,u \rangle=a\}$ be an affine hyperplane ($u\in\R^{n},u\neq 0,
a\in\R$), such that $int(K)\cap H\neq \emptyset$. Then there exist
$z\in int(K)\cap H$ such that $$vol(K)vol(K^{z})\leq
vol(B_{2}^{n})^{2}/4\lambda(1-\lambda),$$ where $\lambda \in[0,1]$
is defined by $\lambda
vol(K)=vol(\{x\in K: \langle x,u\rangle\geq a\})$.
\end{theorem}

The starting point of the results on the {\it Santal\'o-region} of a
convex body is, as stated above, the  {\it Blaschke-Santal\'o inequality}.
This inequality states, roughly, that with the right choice of a
center, the product of the volumes of a body and its polar is
maximized by {\it Euclidean} balls (and also ellipsoids, since this is an
affine invariant). More precisely

\begin{theorem}\label{B}
Let $K$ be a convex body in $\R^{n}$ then \[ \inf _x
Vol(K)Vol(K^{x})\leq Vol(B_{2}^{n})^{2}\] with equality if and only
if $K$ is an ellipsoid.
\end{theorem}

In the special case where $K$ is centrally symmetric, the infimum is attained at $0$, and in 1981 Saint Raymond
\cite{sr} gave a simple proof for this case, namely,   for centrally symmetric $K$ one has
\[ Vol (K) Vol (K^{\circ}) \le Vol (B_2^n)^2,\]
and equality is attained if and only if $K$ is an ellipsoid. The case of equality for a general convex body was proven by Petty \cite{p}.
\begin{rem}\cite{mw}\label{B-infu} For every $K$, a convex body in $\R^{n}$,
there is a unique point $x_{0}$ for which
 $Vol(K)Vol(K^{x_{0}})$ is minimal and this unique $x_{0}$
 is called the {\it Santal\'o point} of $K$.
\end{rem}

We turn now to the {\it Floating body} of convex body.
From  definition, $F(K,\lambda)$ is the
intersection of all half spaces whose  hyperplanes cut off from $K$
a set of volume $\lambda vol(K)$. Note that when $K$ is a centrally
symmetric convex body, we have that its {\it Floating body}
$F(K,\lambda)\not=\emptyset$, for every $0<\lambda<\frac{1}{2}$,
since $0\in F(K,\lambda)$. Moreover, we have $F(K,\lambda)\rightarrow
\{0\}$, in the sense of the {\it Hausdorff distance}, when
$\lambda\rightarrow\frac{1}{2}$ from below. For a general convex
body $K$ we have $F(K,\lambda)$ convergence to $K$ when
$\lambda\rightarrow 0$ from above.

For a centrally symmetric convex body $K$ there is a theorem which was proved by M.Meyer and S.Reisner in
\cite{mr2}. This theorem shows the connection between the supporting hyperplanes of the {\it Floating body}
$F(K,\lambda)$ and the hyperplanes $H_{\lambda,\theta}$, where $\theta\in S^{n-1}$. Before we state this theorem
let us first define $H_{\lambda,K}$ to be the family of all
 hyperplanes $H_{\lambda,\theta}$, where $\theta\in S^{n-1}$.

\begin{theorem}\label{reisner}
Let $K$ be a centrally symmetric convex body in $\R^{n}$, and let $0<\lambda<\frac{1}{2}$. Then there exists a
centrally symmetric convex body $K_{\lambda}$ with the following property: For every supporting hyperplane $H$
of $K_{\lambda}$ we have $vol(K\cap H^{+})= \lambda vol(K)$.
\end{theorem}

For our discussion, we conclude from Theorem \ref{reisner} that the
{\it Floating body} $F(K,\lambda)$, for $K=-K$, is enveloped by the family
$H_{\lambda,K}$. That is, the family of supporting hyperplanes of
$F(K,\lambda)$ coincide with the family $H_{\lambda,K}$, and
$F(K,\lambda)$ lies on the negative side of every $H \in
H_{\lambda,K}$. Theorem \ref{reisner} is a result of a more general
theorem related to measures, namely Theorem \ref{f-is-cvx}.

Before we state Theorem \ref{f-is-cvx}, let us begin with the
following definition related to measures. Let $\mu$ be a
non-negative {\it Borel measure} on $\R^{n}$, and let
$0<\lambda<\frac{1}{2}$. We define $H_{\lambda,\mu}$ to be the
family of all hyperplanes $H$ in $\R^{n}$ such that
$\mu(H^{+})=\lambda$, where as before $H^{+}$ is the half space
bounded by $H$ which does not include $0$.

\begin{theorem}\cite{mr}\label{f-is-cvx}
 Let $\mu$ be a finite, even, log-concave and non-degenerate
measure on $\R^{n}$. Then there exists a convex body $K_{\lambda}$
whose boundary is the envelope of the family $H_{\lambda,\mu }$.
That is, the family of supporting hyperplanes of $K_{\lambda}$
coincides with the family $H_{\lambda,\mu}$ and  $K_{\lambda}$ lies
on the negative side of every $H \in H_{\lambda,\mu}$.
\end{theorem}

\begin{corollary}
Theorem \ref{reisner} follows from Theorem \ref{f-is-cvx} when we
consider  measures of the form $\mu_{K}(L)=vol(K\cap L)$, where $K$
is a centrally symmetric convex body in $\R^{n}$.
\end{corollary}

In order to show the connection between convex bodies and log-concave functions, we need to define the right
analogue for addition and multiplication by scalar of functions. We consider another known operation, which is
related to the {\it Legendre transform}, namely the {\it Asplund product}. Given two functions $f,g:\R^{n}\rightarrow
[0,\infty),$ their {\it Asplund product} is defined by
$$(f*g)(x)=\sup_{x_{1}+x_{2}=x}f(x_{1})g(x_{2}).$$
It is easy to check that $1_{K}*1_{T}=1_{K+T}$ and
$(f*g)^{\circ}=f^{\circ}g^{\circ}$, where $1_{K}$ is the indicator
function. Moreover we argue that the {\it Asplund product} of log-concave
functions is the right analogue for Minkowski addition of convex
bodies.

To define the $\lambda$-homothet of a function $f(x)$, which we denote $(\lambda\cdot f)(x)$, we use
$$(\lambda\cdot f)(x)=f^{\lambda}(\frac{x}{\lambda}).$$ Note that, for a log-concave function, one has indeed
$f*f=2 \cdot f$ and that $(\lambda \cdot f)^{\circ}=(f^{\circ})^{\lambda}.$

To check whether the definitions of Minkowski addition, homothety and duality for log-concave functions are
meaningful and make sense, it is natural to ask whether the basic inequalities for convex bodies such as the
{\it Brunn-Minkowski inequality} and the {\it Santal\'o inequality} remain true, where the role of volume will be played of
course by the integral. We first  discuss the {\it Brunn-Minkowski inequality}.

In its dimension-free form,  the {\it Brunn-Minkowski inequality} states
that $$Vol(\lambda A+ (1-\lambda)B)\geq
Vol(A)^{\lambda}Vol(B)^{1-\lambda},$$ where $A,B\subseteq\R^{n}$
measurable sets and $0<\lambda<1$.

 We have a functional analogue of the {\it Brunn-Minkowski
inequality}, namely the {\it Pr$\acute{e}$kopa-Leindler inequality} (see, e.g. \cite{AKM}) . In the above notation, it
states precisely that
\begin{theorem}
(Pr$\acute{e}$kopa-Leindler) Given $f,g:\R^{n}\rightarrow
[0,\infty)$ and $0<\lambda<1$, $$\int (\lambda \cdot
f)*((1-\lambda)\cdot g)\geq (\int f)^{\lambda}(\int
g)^{1-\lambda}.$$
\end{theorem}

The standard {\it Brunn-Minkowski inequality} follows directly from
{\it Pr$\acute{e}$kopa-Leindler} by considering indicator functions of
sets.

Now, we turn to definition of duality on the class of log-concave functions, its connection with the {\it  Legendre
transform}. As before we  define the dual of a log-concave
function  $f:\R^{n}\rightarrow [0,\infty)$ by
$  f^{\circ}(x)=inf_{y\in \R^{n}}[\frac{e^{-\langle x,y
\rangle}}{f(y)}]$. The above definition of polarity of a log-concave function, is connected with
the {\it Legendre transform} in the following way: $$-log (f^{\circ}) = \L
(-log f),$$ where the {\it Legendre transform} is defined by
$$ (\L \varphi)(x)= \sup_{y \in \R^{n}}(\langle x,y \rangle -\varphi(y)).$$
So, one may conclude that   $$f^{ \circ}= e^{-\L\varphi},$$ where
 $f=e^{-\varphi}$ and $\varphi:\R^{n}\rightarrow\R\cup\{ \infty\}$ is a convex function.

 We define $Cvx(\R^{n})$ to be the class of lower semi-continuous convex
functions $\varphi:\R^{n}\rightarrow\R\cup\{\pm\infty\}$. Let us list some properties of the {\it Legendre transform}.

\begin{enumerate}
\item $\L \varphi \in Cvx(\R^{n})$ for all $\varphi:\R^{n}\rightarrow
\R\cup\{\pm\infty\}$,

\item $\L:Cvx(\R^{n})\rightarrow Cvx(\R^{n})$ is injective and surjective,

\item $\L\L\varphi = \varphi$ for $\varphi\in Cvx(\R^{n})$,

\item $\varphi\leq\psi$ implies $\L\varphi\geq\L\psi$,

\item $\L(max(\varphi,\psi)) = m\hat in(\L\varphi,\L\psi)$ where $ m\hat in$ denotes "regularized minimum",
that is, the largest l.s.c. convex function below all functions
participating in the minimum,

\item $\L\varphi +\L\psi = \L(\varphi \square \psi)$ where $(\varphi  \square \psi)(x):=\inf\{\varphi(y)+\psi(z):x=
y+z\}$.
\end{enumerate}

As the {\it Legendre transform} of any function is a convex function,
$f^{\circ}$ is always log-concave function. In addition if $\varphi$ is
convex and lower semi-continuous then $\L\L \varphi = \varphi$.
Translating to our language, if $f$ is a log-concave upper
semi-continuous function, then $f^{\circ\circ}=f$.

\subsection{Properties of $s$-concave functions}

In what follows, we will use very strongly the notion of an
$s$-concave function, and the fact that such functions can be used
to approximate, quite effectively, log-concave functions. This
technique was used extensively in \cite{AKM} to prove the
{\it Blaschke-Santal\'o inequality} for functions.

We start with the definition of s-concave function. For $s>0$, a function $g:\R^{n}\rightarrow [0,\infty)$ is
called an $s$-concave function if $g^{\frac{1}{s}}$ is concave on the  support of $g$. In other words,
$$\lambda g^{\frac{1}{s}}(x)+(1-\lambda)g^{\frac{1}{s}}(y)\leq g^{\frac{1}{s}}(\lambda x + (1-\lambda)y)
$$ for any $x,y\in supp(g)$ and any $0\leq\lambda\leq1$.

A function on $\R^{n}$ is $s$-concave if and only if it is a
marginal of a uniform measure on a convex body in $\R^{n+s}$.
Indeed, the Brunn concavity principle  implies that such a marginal
is $s$-concave, whereas for the other direction, given an
$s$-concave measure one can easily construct a body $K$ in $R^{n+s}$
which has it as its marginal (there can be many such bodies).

It is easy to see that any $s$-concave function, for $0<s$ integer,
is also log-concave. The converse is of course not true - a gaussian
is log-concave but not $s$-concave for any $s>0$. However, the
family of all functions which are $s$-concave for some $s>0$ is
dense in the family of log-concave functions in any reasonable
sense.

More precisely, let $f:\R^{n}\rightarrow [0,\infty)$ be a
log-concave function. It is possible to create a series, $\{f_{s}
\}_{s=1}^{\infty}$ of $s$-concave functions,  which converges to
$f$, as $s$ tends to infinity, locally uniformly as follows:
$$f_{s}(x):=(1+\frac{logf(x)}{s})^{s}_{+}$$ where $z_{+}=\max\{z,0
\}$. The log-concavity of $f$ implies the s-concavity of $f_{s}$.
Note also that $f_{s}\leq f$ for any $0<s$, and since a log-concave
function is continuous on its support, one has $f_{s}\rightarrow f$
locally uniformly on $\R^{n}$ as $s\rightarrow \infty$.

Let $g:\R^{n}\rightarrow [0,\infty)$ be an $s$-concave function.
Following \cite{AKM}, we define for $s>0$
\begin{eqnarray*}
 \L_{s}g(x)= \inf_{ \{y:g(y)>0  \} } \frac{ (1-
\frac{ \langle x,y \rangle }{s} )^{s}_{+}}{g(y)}.
\end{eqnarray*}
Moreover we define
\begin{eqnarray*}
 K_{s}(g)= \{(x,y)\in \R^{n}\times \R^{s}: \sqrt s
x \in \overline {supp(g)}, |y|\leq g^{\frac{1}{s}}(\sqrt s x) \}.
\end{eqnarray*}
Note that $K_{s}(g)$ is a centrally symmetric convex body in
$\R^{n+s}$ for any even $s$-concave function
$g:\R^{n}\rightarrow[0,\infty)$ such that $supp(g)\subseteq\R^{n}$
is a centrally symmetric convex body.

The following lemma clarifies the relations
between our various definitions.

\begin{lemma}\cite{AKM}\label{C}
For any  $f:\R^{n}\rightarrow [0,\infty),$
\begin{eqnarray*}
K_{s}(f)^{\circ}= K_{s}(\L_{s}(f)) .
\end{eqnarray*}
\end{lemma}

\noindent As appears in \cite{AKM}, it is easy to check that
\begin{itemize}
\item $Vol(K_{s}(g))= \frac{Vol(B_{2}^{s})}{s^{\frac{n}{2}}}  \int
g dx$,
\item   $K_{s}({g_{a}})= K_{s}(g)+
\frac{1}{\sqrt s}(a,0)$.
\end{itemize}

 For  $s$-concave function $g:\R^{n}\rightarrow [0,\infty)$, we call $\L_{s} g$ the polar of
$g$. From the definition of $\L_{s}$ it is easy to see that $\L_{s}g\leq g^{\circ}$,  and $\L_{s}\L_{s} g =g$
for any upper semi-continuous, $s$-concave function $g:\R^{n}\rightarrow [0,\infty)$ which satisfies $g(0)>0$ (see \cite{am}).

As may seem natural by now, we have a {\it Santal\'{o}-type inequality}
also for $s$-concave functions, where we use $\L_{s}$ for the
duality operation. This inequality
states that

\begin{theorem}\cite{AKM}\label{k1}
Let $g$ be an $s$-concave function on $\R^{n}$, with $0<\int
g<\infty$, whose center of mass is at the origin (i.e., $\int
xg(x)=0$). Then
$$\int_{\R^{n}}g\int_{\R^{n}}\L_{s}(g)\leq \frac{s^{n}Vol(B_{2}^{n+s})^{2}}{Vol(B_{2}^{s})^{2}},$$
with equality if and only if $g$ is a marginal of a uniform
distribution of an $(n+s)$-dimensional ellipsoid. And
$B_{2}^{n+s}$,$B_{2}^{s}$ are the {\it Euclidean} balls of radius 1 in
$\R^{n+s} $ and $\R^{s}$ respectively.
\end{theorem}

\newpage

\section{The Santal\'o  Region of a log-concave function}\label{TSR}

Let $\varphi:\R^{n}\rightarrow
(-\infty,\infty]$ be a  convex function.   We discuss some properties of the set
$$ S(e^{-\varphi},t)=\bigg\{a  \in  \R^{n}  :  \int e^{-\varphi_{a}} \int e^{-\L \varphi_a} \le  (2\pi)^n t \bigg\}.$$
 According to \cite{mw}, we call this set the
{\it Santal\'{o} region} of the function $e^{-\varphi}$, with constant
$t$. It follows from the functional version of the {\it Santal\'{o}
inequality}, that the set $S(e^{-\varphi}, 1)$ is non-empty for every
convex $\varphi$, see \cite{AKM}.

 Note that the {\it Legendre
transform} of the translate $\varphi_a$ is simply
\begin{eqnarray*} (\L \varphi_a) (x) &=&
\sup_{y \in \R^n} \left( \langle x, y \rangle - \varphi_a(y) \right)
=
\sup_{y \in \R^n} \left( \langle x, y \rangle - \varphi (y -a ) \right) \\
 &= & \sup_{z \in \R^n} \left( \langle x, z+a \rangle - \varphi(z) \right)
 = \langle x, a \rangle + \sup_{z \in \R^n} \left( \langle x, z
\rangle - \varphi(z) \right)   \\ & = &  \langle x, a \rangle  + (\L
\varphi)(x).
\end{eqnarray*}

\noindent Therefore,  using that $\int e^{-\varphi} = \int e^{-\varphi_a}$ for
every $a$, we may write
 \[ S(e^{-\varphi},t)=\bigg\{a \in
\R^{n} : \int e^{-\langle a,x\rangle }e^{-(\L \varphi)(x)} dx \leq
\frac{(2 \pi)^{n}t}{ \int e^{-\varphi}dx }  \bigg\}. \]
For every $0<\lambda<1$,  $a,b\in \R^{n}$
\[
\begin{array}{l}
(\L \varphi_{\lambda a + (1-\lambda) b })(x) = \langle x,\lambda a
+(1-\lambda) b)\rangle + (\L \varphi)(x)
\\
\\
= \lambda (\langle x,a
\rangle + (\L \varphi)(x)) + (1- \lambda)(\langle x,b\rangle + (\L
\varphi)(x))
\\
\\
=\lambda (\L \varphi_{a}(x)) + (1-\lambda)(\L \varphi_{b}(x)).
\end{array}
\]

If $a \neq b$ in $\R^{n}$ and $\varphi \not\equiv 0$ (meaning there
exists $x\in \R^{n}$ such that $\varphi(x)\neq 0$) then there exist
$x_{0}\in \R^{n}$ such that $\L \varphi_{a}(x_{0}) \neq \L
\varphi_{b}(x_{0})$. Indeed , we have
\begin{eqnarray*}
\L \varphi_{a}(x_{0}) = \langle x_{0}, a \rangle + (\L \varphi)(x_{0}) \\
\L \varphi_{b}(x_{0}) = \langle x_{0}, b \rangle + (\L
\varphi)(x_{0}).
\end{eqnarray*}
So one may choose  $x_{0}\in \R^{n}$ such that $\langle
x_{0},a \rangle \neq \langle x_{0},b \rangle$ and $\L
\varphi(x_{0})\neq \infty$.

Define $F: \R^n \to \R^+$ by $F(a)= \int e^{- \L \varphi_{a}(x)}dx$.
Note that $F$ is  a strictly convex function. Indeed
\begin{eqnarray*}
 \int e^{-\L \varphi_{\lambda a +( 1 - \lambda )b}(x)}dx = \int
(e^{ -\L \varphi_{a}(x)})^ \lambda (e^{- \L \varphi_{b}(x)})^{1-
\lambda}
\\< \lambda \int e^{-\L \varphi_{a}(x)}dx + ( 1 - \lambda ) \int e^{-\L
\varphi_{b}(x)}dx,
\end{eqnarray*}
where the last inequality follows from the {\it arithmetic-geometric mean
inequality} and the existence of $x_{0} \in \R^{n}$ such that $\L
\varphi_{a}(x_{0}) \neq \L \varphi_{b}(x_{0})$, where $a
\neq b$ and $\varphi\not\equiv0$.

Consider $ \varphi\not\equiv0$.  We call a point $a\in \R^n$ a
{\it Santal\'{o} point } of $e^{-\varphi}$ if $F(a) \leq F(b)$ for every $b$. Since $F$ is a strictly convex function, the {\it Santal\'o point } is
unique. It is easy to see that $F$ is a continuous function on
$\{x\in\R^{n}:F(x)<\infty\}$ (Lemma 3.2 in \cite{AKM}) and differentiable (since the differentiation under the
integral sign is allowed, as the derivative is locally bounded and integrals converge), and that $\nabla
F(a)=\int xe^{- \L \varphi_{a}(x)}dx$, which is equal to $0$ if and only if $a$ is the {\it Santal\'o point} of
$e^{-\varphi}$. Moreover, one may conclude that $\grad F(b)\neq 0$ for every $b\in\partial(S(e^{-\varphi},t))$
whenever $int(S(e^{-\varphi},t))\neq \emptyset$. From  strict convexity of $F$ we have  $0
\in int(S(f,t))$ whenever $int(S(f,t))\neq {\emptyset}$.

Unless stated otherwise we will always assume that $e^{-\varphi}$
has its {\it Santal\'{o}-point} at the origin and that $\int
e^{-\varphi}dx < \infty$.

\begin{lemma}\label{1}
For every convex function $\varphi: \R^{n}\rightarrow \R \cup \{
+\infty \} $ and any $t> 0$, the {\it Santal\'{o} region}
 $S(e^{-\varphi},t)$ is either a strictly convex body, or the empty set.
\end{lemma}

\noindent {\bf Proof.} Consider  $t>0$. Let $a\neq b$ and for them to be two points at the
boundary of $S(e^{-\varphi}, t)$. We have to show that for any
$0<\lambda<1$, the point $\lambda a + (1-\lambda)b$ is in the
interior of $S(e^{-\varphi}, t)$. That is, we know that
$F(a)=F(b)=\frac{(2\pi)^{n}t}{\int e^{-\varphi}}$ and need to show
that $F(\lambda a +(1-\lambda)b)<\frac{(2\pi)^{n}t}{\int
e^{-\varphi}}$, which follows from strict convexity of $F$.\qed

\begin{Fact}
 The mapping $t \rightarrow S(e^{-\varphi},t) $ is increasing and it is also concave,
 i.e for every $0<\lambda<1$ we have $ \lambda S(e^{-\varphi},t)+ (1-\lambda)S(e^{-\varphi},s)\subseteq S(e^{-\varphi},\lambda t
 +(1-\lambda)s)$.
\end{Fact}

\noindent {\bf Proof.} The fact that the map is increasing follows
immediately from  definition. To prove the concavity consider $a
\in S(e^{-\varphi},t)$ and $ b \in S(e^{-\varphi},s) $. We have
\begin{eqnarray*}
F(a) \leq \frac{(2 \pi)^{n}t}{\int e^{-\varphi}}dx \ , \ F(b) \leq
\frac{(2 \pi)^{n}s}{\int e^{-\varphi}}dx.
\end{eqnarray*}
Thus, from strict convexity of $F$ we conclude
$$F(\lambda a + (1-\lambda)b)< \lambda F(a) + (1-\lambda)F(b)\leq (\lambda t + (1-\lambda)s) \frac{(2\pi)^{n}}{\int e^{-\varphi}}. $$ \qed

Let us check how linear transformations affect  the {\it Santal\'o
region}. Let $\varphi:\R^{n}\rightarrow \R\cup \{\infty\}$ be a
convex function and let $A$ be an invertible linear map, $A\in
GL_{n}$.

From definition
\[
\begin{array}{l}
 A(S(e^{-\varphi},\lambda))= \{a\in\R^{n}:\int e^{-\varphi(x)}dx \int
e^{-\L \varphi_{A^{-1}a}(x)}dx \leq (2\pi)^{n}\lambda \},
\\
\\
 S(e^{-\varphi\circ A^{-1}},\lambda)= \{a\in\R^{n}: \int e^{-\varphi\circ A^{-1}(x)}dx
\int e^{-\L ((\varphi\circ A^{-1})_{a})(x)}dx\leq (2\pi)^{n}\lambda
\}.
\end{array}
\]
Note that these two regions  coincide.

\begin{Fact}\label{lemma3}
For every convex function $\varphi: \R^{n} \rightarrow \R \cup \{
+\infty \}$,   the boundary of $S(e^{-\varphi},t)$, denoted by $
\partial (S(e^{-\varphi},t))$, is infinitely  smooth, whenever $int (S(e^{-\varphi},t))\neq\emptyset$.
\end{Fact}

\noindent {\bf Proof.} Assume that  $int
(S(e^{-\varphi},t))\neq\emptyset$, and define as before $F(a)= \int
e^{-\L \varphi_{a}(x)dx}$ . The boundary of $S(e^{-\varphi},t)$ is
given by
\begin{eqnarray*}
\bigg\{ a \in \R^{n} : F(a) = \frac{(2\pi)^{n}t}{\int e^{-\varphi}dx} \bigg\}
=
\partial (S(e^{-\varphi},t)).
\end{eqnarray*}
Differentiation under the integral sign gives
\begin{eqnarray*}
\frac{ \partial F(a) }{\partial a_{j} } = (-1) \int x_{j} e^{\langle
x,a \rangle}e^{- \L \varphi(x)}dx \textrm{ for every } 1\leq j \leq n.
\end{eqnarray*}
Note that the differentiation under
the integral sign is allowed, as the derivative is locally bounded
and the integrals converge. Thus as we said before we have  $\nabla
F(b) \neq 0 $ for every $b\in
\partial (S(e^{-\varphi},t)) $.
Therefore Fact \ref{lemma3} follows from the {\it implicit function
theorem}.\qed

\newpage

\section{The relation between the Santal\'o region and
the Floating body of an even log-concave function}
\label{bigsection}

Throughout this section we consider only  even log-concave functions. We will compare two different bodies, the
{\it Santal\'o region} and the {\it Floating body}, related to a specific even log-concave function. For any positive integer $s$ we consider $\R^n \subset \R^{n+s}$ as the first $n$ coordinates. Moreover,  we denote $P_{\R^n}$ as the projection onto $\R^n$.

We start with some basic notations. Let $f=e^{-\varphi}$ be an even log-concave function and let $t>0$. We define
$$A_{t}=\{x\in \R^{n}:\varphi(x)\leq t \}.$$ It is easy to check that $A_{t}$ is a centrally symmetric body
which is convex in $\R^{n}$, also since we require  $\int f < \infty $, we get that $A_{t}$ is a bounded
set. We define

 \begin{displaymath} \varphi_{[t]}(x) := \left  \{
\begin{array}{ll}
\varphi(x)  & \textrm{if $x\in A_{t}$}\\
\infty & \textrm{if $x \not\in A_{t}$},\\
\end{array} \right.
\end{displaymath}

 $$f_{[t]}(x):=e^{-\varphi_{[t]}(x)}.$$

\noindent Obviously $f_{[t]}$ is an even log-concave function such that
$supp(f_{[t]})=A_{t}$. In addition, we have  $f_{[t]}(x)\geq{e^{-t}}>0$
for any $x\in supp(f_{[t]})$.

\noindent As before we define
$f_{[t],s}(x)=(1+\frac{logf_{[t]}(x)}{s})^{s}_{+}$. Note that
\begin{itemize}
\item $supp(f_{[t],s})=supp(f_{[t]})$,

\item  $(f_{[t]})_{s}(x)>\frac{e^{-t}}{2}>0$ for every $x\in
supp(f_{[t],s})$  for $s$ big enough.

\end{itemize}
In addition,  for $a\in \R^{n}$ we
define $$f_{[t],s,a}(x):= f_{[t],s}(x-a).$$

It is obvious that for any even log-concave function $f$, $f_{[t]}$ is an even log-concave function for any positive $t$, with compact support, which converges locally uniformly to $f$ as $t\rightarrow\infty$.  Our goal is to prove Theorem \ref{Thm} and Theorem \ref{reverse} for any  function of the form  $f_{[t]}$,  $t>0$. Then from continuity (with respect to {\it Hausdorff distance} of bodies) of the {\it Floating body} and the {\it Santal\'o region}  with respect  to locally uniformly convergence of even log-concave functions (Lemma \ref{c1} and Lemma \ref{c2}) we will complete the proof of Theorem \ref{Thm} and Theorem \ref{reverse} for any even log-concave function.

The following lemma will be the main ingredient in the proof of Theorem \ref{Thm} and Theorem \ref{reverse}.
As mentioned above,  Lemma \ref{c1} deals with continuity of  the {\it Floating body}   with respect  to locally uniformly convergence of even log-concave functions.

\begin{lemma} \label{c1}
Let $f:\R^{n}\rightarrow [0,\infty)$ be an even log-concave
function. Then
$$F(f_{[t]},\lambda)\rightarrow_{t\rightarrow\infty}F(f,\lambda) \textrm{ for any } 0<\lambda<1/2.$$
\end{lemma}

\noindent {\bf Proof.} Let $\theta\in S^{n-1}$ and let $H_{f}$, $H_{f_{[t]}}$ be two hyperplanes  perpendicular
to $\theta$ and such that $\int_{H_{f}^{+}} f dx = \lambda \int f dx$, and
 $\int_{H_{f_{[t]}}^{+}} f_{[t]} dx = \lambda \int f_{[t]} dx$.

From definition of $f_{[t]}$, it is easy to see that $\lambda \int f_{[t]}dx \rightarrow_{t\rightarrow\infty}\lambda \int f dx$. Since
 $H_{f}$ and  $H_{f_{[t]}}$ are parallel, we conclude
$H_{f_{[t]}}\rightarrow_{t\rightarrow\infty}H_{f}$ with respect to the {\it Hausdorff distance} of bodies, which
implies  $F(f_{[t]},\lambda)\rightarrow_{t\rightarrow\infty}F(f,\lambda)$.\qed

Lemma \ref{c2} deals with the continuity of the  {\it Santal\'o region}   with respect  to locally uniformly convergence of even log-concave functions.

\begin{lemma}\label{c2}
Let $f:\R^{n}\rightarrow [0,\infty)$ be an even log-concave
function. Then
$$S(f_{[t]},d)\rightarrow_{t\rightarrow\infty}S(f,d)\textrm { for any } d>0.$$
\end{lemma}

\noindent {\bf Proof.}  Consider, as before,  $f=e^{ -\varphi}$, where $\varphi$ in an even convex function, and $f_{[t]}= e^{ -\varphi_{[t]}}$.  From definition we know
$$
 S(f_{[t]},d)= \big\{a \in \R^{n}: \int f_{[t]}  \int
 ((f_{[t]})_{a})^{\circ} \leq (2 \pi)^{n}d \big\},$$
$$S(f,d)= \big\{a \in \R^{n}: \int f \int
 (f_{a})^{\circ} \leq (2 \pi)^{n}d \big\}.$$

\noindent Moreover, we define
\begin{itemize}
\item $F_{t}(a):=\int
 ((f_{[t]})_{a})^{\circ}dx=\int e^{-\langle x,a \rangle}e^{-\L (\varphi_{[t]})(x)}dx$,
\item  $F(a):=\int
 (f_{a})^{\circ}dx=\int e^{-\langle x,a \rangle}e^{-\L (\varphi)(x)}dx$.
\end{itemize}

Our general idea is as follows. First  we  prove $F_{t}\rightarrow_{t\rightarrow\infty}F$ locally
 uniformly on $\{x\in\R^{n}: F(x)<\infty \}$. Then, since $S(f,d)\subseteq \{x\in\R^{n}: F(x)<\infty
 \}$ is a closed bounded set, from strict convexity and continuity  of $F$
 in $S(f,d)$ (Lemma 3.2 \cite{AKM})  we  complete the proof.

In order to show  locally uniformly  convergence, let us first prove
 $F_{t}(a)\rightarrow_{t\rightarrow\infty}F(a)$ for any $a\in \{x\in\R^{n}: F(x)<\infty
 \}$. It is easy to check that $\L (\varphi_{[t]})(x)\rightarrow_{t\rightarrow\infty}\L
 (\varphi)(x)$ for any $x\in\R^{n}$, hence  $$e^{-\langle x,a \rangle}e^{-\L (\varphi_{[t]})(x)} \rightarrow_{t\rightarrow\infty}e^{-\langle x,a \rangle}e^{-\L (\varphi)(x)}\textrm{ for every  } x\in \R^{n}.$$

\noindent Since $e^{-\langle x,a \rangle}e^{-\L
(\varphi_{[t]})(x)}$ and $e^{-\langle x,a \rangle}e^{-\L
(\varphi)(x)}$ are both log-concave functions (Lemma 3.2 \cite{AKM}),  we conclude
$F_{t}(a)\rightarrow_{t\rightarrow\infty}F(a)$.

Moreover, for any $a\in\{x\in\R^{n}:F(x)<\infty \}$ it is easy to see that $F_{t_{2}}(a) \leq F_{t_{1}}(a)$ whenever
$0<t_{1}<t_{2}$.   Indeed, since  $
\varphi_{[t_{2}]}(x)\leq\varphi_{[t_{1}]}(x)$ for any $x\in \R^{n}$ one has
 $\L (\varphi_{[t_{1}]})(x)\leq
\L(\varphi_{[t_{2}]})(x)$ for any $x\in \R^{n}$.

 Since   $F_{t}\rightarrow_{t\rightarrow\infty} F$ pointwise in $\{x\in\R^{n}: F(x)<\infty \}$,    and  $(F_{t})_{t>0}$ is a decreasing sequence of functions, we conclude from  continuity of $F$ in $\{x\in\R^{n}: F(x)<\infty \}$ that
$F_{t}\rightarrow_{t\rightarrow\infty}F$ locally
 uniformly in $\{x\in\R^{n}: F(x)<\infty \}$. In addition, since  $F$ is a strictly convex function in $S(f,d)$, then $S(f,d)\subseteq\{x\in\R^{n}: F(x)<\infty
 \}$ is a bounded
convex set in $\R^{n}$. \qed

Let us prove the following lemma.
\begin{lemma}\label{10}
Let $f:\R^{n}\rightarrow [0,\infty)$ be an even log-concave function
such that $supp(f)$ is a convex compact body in $\R^{n}$, and there
exist $0<t_{0}$ such that $ t_{0}<f(x)$ for all $x\in supp(f)$.
Define as before $f_{s}(x)=(1+ \frac{logf(x)}{s})_{+}^{s}$ and let
$0<d$. Then
\begin{eqnarray*}
 \bigg\{a \in \R^{n}: \int f_{s}  \int
 \L_{s}((f_{s})_{a}) \leq (2 \pi)^{n}d \bigg\} \rightarrow_{s\rightarrow \infty}
 S(f,d),
\end{eqnarray*}
with respect to the {\it Hausdorff distance} of convex bodies.
\end{lemma}

\begin{corollary}
From Lemma \ref{10} we conclude that $$\bigg\{a \in \R^{n}: \int
f_{[t],s} \int
 \L_{s}(f_{[t],s,a}) \leq (2 \pi)^{n}d \bigg\} \rightarrow_{s\rightarrow \infty}
 S(f_{[t]},d)$$
for any even log-concave function $f:\R^{n}\rightarrow[0,\infty)$
and any $0<t$.
\end{corollary}

 \noindent {\bf Proof of Lemma \ref{10}.} For $s$ big enough, $supp(f_{s})=supp(f)$. We define
\begin{itemize}
\item $F(a):= \int_{\R^{n}} (f_{a})^{\circ}(x)dx= \int_{\R^{n}} \inf_{ y\in supp(f)+a}\frac{e^{-\langle x,y \rangle}}{f(y-a)}dx$,
\item $F_{s}(a) := \int_{\R^{n}} \L_{s}((f_{s})_{a}) dx= \int_{\R^{n}} \inf_{y\in supp(f)+a}\frac{(1-\frac{\langle x,y\rangle}{s})_{+}^{s}}{f_{s}(y-a)}dx$.
\end{itemize}

First, we prove  $F_{s}\rightarrow_{s\rightarrow\infty}F$ locally uniformly in $\{x\in\R^{n}:F(x)<\infty \}$. Then, since $S(f,d)\subseteq\{x\in\R^{n}:F(x)<\infty \}$ is a bounded convex set in $\R^{n}$,  from the fact that
$F$ is  a continuous and  strictly convex function in $S(f,d)$, we complete the proof.

Let $D\subseteq\{x\in\R^{n}:F(x)<\infty \}$ be a closed, bounded set
and let $\epsilon>0$. We will prove that there exist
$s_{\epsilon}>0$ such that for every $s\ge s_{\epsilon}$
$$|F_{s}(a)-F(a)|<\epsilon \textrm{ for any }a\in D.$$

Since  $F$ is a continuous function in  $
\{x\in\R^{n}:F(x)<\infty \} $, there exist
$r_{\epsilon,D}>0$ such that
\begin{eqnarray} \label{a1}
 |\int_{\R^{n}-r_{\epsilon,D}B_{2}^{n}} \inf_{ y\in supp(f)+a}\frac{e^{-\langle x,y
\rangle}}{f(y-a)}dx |< \frac{\epsilon}{8}, \textrm{ for any } a\in D.
\end{eqnarray}

\noindent For this moment let's assume  there exists $s_{\epsilon}>0$ such that for
every $s\ge s_{\epsilon} $
\begin{eqnarray} \label{b1}
|\int_{r_{\epsilon,D}B_{2}^{n}} \inf_{ y\in supp(f)+a}\frac{e^{-\langle x,y
\rangle}}{f(y-a)}dx - \int_{r_{\epsilon,D}B_{2}^{n}} \inf_{y\in
supp(f)+a}\frac{(1-\frac{\langle
x,y\rangle}{s})_{+}^{s}}{f_{s}(y-a)}dx| < \frac{\epsilon}{4}
\end{eqnarray}
and
\begin{eqnarray} \label{cc1}
|\int_{\R^{n}- r_{\epsilon,D}B_{2}^{n}} \inf_{y\in
supp(f)+a}\frac{(1-\frac{\langle
x,y\rangle}{s})_{+}^{s}}{f_{s}(y-a)}dx| <\frac{\epsilon}{4}
\textrm{ for any } a\in D.
\end{eqnarray}
 From (\ref{a1}),(\ref{b1}),(\ref{cc1}) and  the {\it triangle-inequality}  we get
\[
\begin{array}{l}|\int_{\R^{n}} \inf_{ y\in
supp(f)+a}\frac{e^{-\langle x,y \rangle}}{f(y-a)}dx - \int_{\R^{n}}
\inf_{y\in supp(f)+a}\frac{(1-\frac{\langle
x,y\rangle}{s})_{+}^{s}}{f_{s}(y-a)}dx| \leq
\\
\\
\leq |\int_{r_{\epsilon,D}B_{2}^{n}} \inf_{ y\in
supp(f)+a}\frac{e^{-\langle x,y \rangle}}{f(y-a)}dx -
\int_{r_{\epsilon,D}B_{2}^{n}} \inf_{y\in
supp(f)+a}\frac{(1-\frac{\langle
x,y\rangle}{s})_{+}^{s}}{f_{s}(y-a)}dx|
\\
\\
+|\int_{\R^{n}-r_{\epsilon,D}B_{2}^{n}} \inf_{ y\in
supp(f)+a}\frac{e^{-\langle x,y \rangle}}{f(y-a)}dx | +
|\int_{\R^{n}- r_{\epsilon,D}B_{2}^{n}} \inf_{y\in
supp(f)+a}\frac{(1-\frac{\langle
x,y\rangle}{s})_{+}^{s}}{f_{s}(y-a)}dx|< \epsilon,
\end{array}
\]
which shows $F_{s}\rightarrow_{s\rightarrow\infty}F$ locally uniformly in
$\{x\in\R^{n}:F(x)<\infty \}$.

Let us show  (\ref{b1}). Let $x\in r_{\epsilon,D}B_{2}^{n}$ , $a\in D$ and $y\in supp(f)+a\subseteq
supp(f)+D$. Since $t_{0}<f$ we get
$$|\frac{e^{-\langle x,y
\rangle}}{f(y-a)}- \frac{(1-\frac{\langle x,y\rangle}{s})_{+}^{s}}{f_{s}(y-a)}|\leq  \frac{1}{t_{0}}
|e^{-\langle x,y \rangle}- \frac{f(y-a)}{f_{s}(y-a)} (1-\frac{\langle x,y\rangle}{s})_{+}^{s}|.$$

 Now, since
$r_{\epsilon,D}B_{2}^{n}$ and $supp(f)+D$ are bounded sets,  the value of $\langle x,y \rangle$ is bounded whenever $x\in r_{\epsilon,D}B_{2}^{n}$ and $y\in
supp(f)+D$. Since $f_{s}\rightarrow_{s\rightarrow\infty}f$ locally uniformly in $supp(f)$, and since $\frac{t_{0}}{2}<f_{s}$ for
$s$ is big enough, we conclude that for any $\epsilon_{0}>0$ there exists $s_{\epsilon_{0}}>0$ such that for any $
s\ge s_{\epsilon_{0}}$
$$\frac{1}{t_{0}}|e^{-\langle x,y \rangle}-
\frac{f(y-a)}{f_{s}(y-a)} (1-\frac{\langle x,y\rangle}{s})_{+}^{s}|<
\epsilon_{0} $$ for any $x\in r_{\epsilon,D}B_{2}^{n}$ , $a\in D$
and $y\in supp(f)+a.$

Thus $$ |\inf_{y\in
supp(f)+a}\frac{e^{-\langle x,y \rangle}}{f(y-a)}- \inf_{y\in
supp(f)+a}\frac{(1-\frac{\langle
x,y\rangle}{s})_{+}^{s}}{f_{s}(y-a)}|\leq \epsilon_{0}.$$
Choose $0<\epsilon_{0}< \frac{\epsilon}{4
r_{\epsilon,D}^{n}Vol(B_{2}^{n})}$. Then
\[
\begin{array}{l}|\int_{r_{\epsilon,D}B_{2}^{n}} \inf_{ y\in
supp(f)+a}\frac{e^{-\langle x,y \rangle}}{f(y-a)}dx -
\int_{r_{\epsilon,D}B_{2}^{n}} \inf_{y\in
supp(f)+a}\frac{(1-\frac{\langle
x,y\rangle}{s})_{+}^{s}}{f_{s}(y-a)}dx|
\\
\\
\leq\int_{r_{\epsilon,D}B_{2}^{n}} | \inf_{ y\in
supp(f)+a}\frac{e^{-\langle x,y \rangle}}{f(y-a)}-\inf_{y\in
supp(f)+a}\frac{(1-\frac{\langle
x,y\rangle}{s})_{+}^{s}}{f_{s}(y-a)}|dx
\\
\\
\leq r_{\epsilon,D}^{n}Vol(B_{2}^{n}) \epsilon_{0}< \frac{\epsilon}{4}
\end{array}
\]
We continue with (\ref{cc1}). We know $(1-\frac{\langle
x,y\rangle}{s})_{+}^{s}\leq e^{-\langle x,y \rangle}$, which implies
\[
\begin{array}{l}
|\int_{\R^{n}- r_{\epsilon,D}B_{2}^{n}} \inf_{y\in
supp(f)+a}\frac{(1-\frac{\langle
x,y\rangle}{s})_{+}^{s}}{f_{s}(y-a)}dx|
\\
\\
\leq  |\int_{\R^{n}-
r_{\epsilon,D}B_{2}^{n}} \inf_{y\in supp(f)+a}\frac{e^{-\langle
x,y\rangle}}{f_{s}(y-a)}dx|.
\end{array}
\]

Since $f=e^{-\varphi}$ and
$f_{s}=e^{-\varphi_{s}}$ where
$\varphi_{s},\varphi:\R^{n}\rightarrow \R\cup\{\infty\}$ are convex
functions, and since $f_{s}\leq f$ and $f_{s}\rightarrow f$ locally
uniformly, we conclude that $\varphi\leq\varphi_{s}$ and that
$\varphi_{s}\rightarrow \varphi$ locally uniformly in $supp(f)$.
Since $supp(f)$ is a bounded set in $\R^{n}$,
for $\epsilon_{1}>0$  there   exists $s_{\epsilon_{1}}>0$ such that for
any $s\ge s_{\epsilon_{1}}$
 $$ \varphi_{s}(x)\leq
\varphi(x)+\epsilon_{1}\textrm{ for every } x\in supp(f).$$ Hence
 $$-\L( \varphi_{s})(x)\leq -\L (\varphi+\epsilon_{1})(x) =
\epsilon_{1} - \L(\varphi)(x) \textrm{ for every } x\in supp(f).$$

\noindent  From  (\ref{a1})
\[
\begin{array}{l}
\int_{\R^{n}-
r_{\epsilon,D}B_{2}^{n}} \inf_{y\in supp(f)+a}\frac{e^{-\langle x,y
\rangle}}{f_{s}(y-a)}dx
\\
\\
 = \int_{\R^{n}- r_{\epsilon,D}B_{2}^{n}}
e^{-\langle x,a \rangle}
 e^{-\L (\varphi_{s})(x)}dx
\\
\\
\leq  e^{\epsilon_{1}}\int_{\R^{n}- r_{\epsilon,D}B_{2}^{n}} e^{-\langle x,a \rangle}
e^{-\L (\varphi)(x)}dx  < e^{\epsilon_{1}} \frac{\epsilon}{8}.
\end{array}
\]
We choose $0<\epsilon_{1}<ln(2)$ and  $s_{\epsilon}= \max
\{s_{\epsilon_{0}},s_{\epsilon_{1}} \}$.
\qed

As mentioned before,  any even $s$-concave function $f_s:\R^{n}\rightarrow [0,\infty)$
such that $supp(f_{s})$ is convex, bounded and symmetric with respect to the origin, one has
$K_{s}(f_{s})\subseteq\R^{n+s} $ is a centrally symmetric convex body. Lemma \ref{projection} shows the
connection between the projection of $F (K_{s}(f_{s}),\lambda)$ on $\R^{n}$ ($\R^{n}\subseteq \R^{n+s}$ the
first $n$ coordinates) and  the projection of the intersection   of all half spaces which contain $F
(K_{s}(f_{s}),\lambda)$  and also  determined by a hyperplane in $\R^{n+s}$ that is perpendicular to $\R^{n}$
and supports $F (K_{s}(f_{s}),\lambda)$.

\begin{lemma}\label{projection}
Let $f:\R^{n}\rightarrow [0,\infty)$ be an even log-concave
function. Define as before  $f_{s}(x)=(1+
\frac{logf(x)}{s})_{+}^{s}$. Then
\begin{eqnarray*}
P_{\R^{n}}[ F (K_{s}(f_{s}),\lambda)]=
P_{\R^{n}}[\bigcap_{\theta=(\theta_{1},..,\theta_{n},0,..,0) \in
S^{n+s-1} } H_{\lambda,\theta}^{-}],
\end{eqnarray*}
where $\lambda
Vol(K_{s}(f_{s}))=Vol(K_{s}(f_{s})\cap H_{\lambda,\theta}^{+})$   for any $\theta =
(\theta_{1},..,\theta_{n},0,..,0 ) \in S^{n+s-1}$.
\end{lemma}

\noindent {\bf Proof.} Let $H$ be an affine hyperplane in $\R^{n+s}$ ($1\leq dim(H)\leq n+s-1$) and let $K$ be a
centrally symmetric convex body in $\R^{n+s}$. We denote by $P_{H}(K)$ as the projection of $K$ on $H$. It is
easy to check
\begin{eqnarray}
P_{H}(K)=P_{H}(\bigcap_{\{\theta\in S^{n+s-1}\cap H   \}}
H_{\theta}^{-}),
\end{eqnarray}
where $H_{\theta}$ is a supporting hyperplane to $K$ and
perpendicular to $\theta$ such that $K\subseteq H_{\theta}^{-}$.

From Theorem \ref{f-is-cvx},  $F (K_{s}(f_{s}),\lambda)$
is  enveloped by all hyperplanes $H_{\lambda,\theta}$ such that
\begin{eqnarray*}
Vol_{n+s}(K_{s}(f_{s}) \cap H_{\lambda,\theta}^{+}) = \lambda
Vol_{n+s}(K_{s}(f_{s})) .
\end{eqnarray*}
Meaning, any hyperplane $\widetilde{H}$, such that $F
(K_{s}(f_{s}),\lambda)\subseteq \widetilde{H}^{-}$, is a supporting
hyperplane to $F(K_{s}(f_{s}),\lambda)$ if and only if
$\widetilde{H}$ satisfies
$Vol( K_{s}(f_{s})\cap \widetilde{H}^{+} )= \lambda
Vol(K_{s}(f_{s}))$. We take $K=F(K_{s}(f_{s}),\lambda)$, $H=\R^{n}$
($\R^{n}\subset\R^{n+s}$ the first $n$ coordinates). \qed

The next lemma shows the connection between the {\it Floating body} which
is related to a convex body  and the {\it Floating body} which is related
to a $s$-concave function.

\begin{lemma}\label{good}
Let $f:\R^{n}\rightarrow [0,\infty)$ be an even log-concave function
and let $0<\lambda<\frac{1}{2}$. Define as before  $f_{s}(x)=(1+
\frac{logf(x)}{s})_{+}^{s}$. Then
 $$F(f_{s},\lambda)= \sqrt s P_{\R^{n}}[  F (K_{s}(f_{s}),\lambda)].$$

\end{lemma}

\noindent {\bf Proof.} Let $0<\lambda<\frac{1}{2}$ and $  \theta
=(\theta_{1},..,\theta_{n},0,..,0) \in S^{n+s-1}  $. We will compute
$Vol(K_{s}(f_{s})\cap H_{\lambda,\theta}^{+})$ in two different
ways. On the one hand
\begin{eqnarray*}
Vol(K_{s}(f_{s})\cap H_{\lambda,\theta}^{+}) = \lambda
Vol(K_{s}(f_{s}))= \lambda \frac{Vol(B_{2}^{s})}{s^{\frac {n}{2}}}
\int f_{s}dx.
\end{eqnarray*}
On the other hand, simple change of variables $x\rightarrow \frac{1}{\sqrt s}x$ yields
\[
\begin{array}{l}
Vol(K_{s}(f_{s})\cap H_{\lambda,\theta}^{+})
=\int_{P_{\R^{n}}(H_{\lambda,\theta})^{+} } Vol
(f_{s}^{\frac{1}{s}}(\sqrt s x)B_{2}^{s})dx
\\
\\
=\frac{Vol(B_{2}^{s})}{s^{\frac{n}{2}}}\int _{\sqrt s
P_{\R^{n}}(H_{\lambda,\theta})^{+}} f_{s}(x)dx.
\end{array}
\]
 We conclude
\begin{eqnarray*}
\lambda  \int f_{s}dx=\int _{\sqrt s
P_{\R^{n}}(H_{\lambda,\theta})^{+}} f_{s}(x)dx, \  \ \forall  \theta = (\theta_{1},..,\theta_{n},0,..,0) \in S^{n+s-1}.
\end{eqnarray*}

\noindent From  definition of  {\it Floating body} and  Lemma \ref{projection}  we get
\[
\begin{array}{l}
F(f_{s},\lambda)= \bigcap
_{\theta=(\theta_{1},..,\theta_{n},0..,0)}\sqrt s
P_{\R^{n}}(H_{\lambda,\theta})^{-}
\\
\\
=  \sqrt s \bigcap
_{\theta=(\theta_{1},..,\theta_{n},0..,0)}
P_{\R^{n}}(H_{\lambda,\theta}^{-})
\\
\\
= \sqrt s
P_{\R^{n}}[\bigcap_{\theta=(\theta_{1},..,\theta_{n},0..,0)}H_{\lambda,\theta}^{-}]
\\
\\
=\sqrt s P_{\R^{n}}[F(K_{s}(f_{s}),\lambda)].
\end{array}
\]
 \qed \\

Let $K\subseteq \R^{n+s}$ be a centrally symmetric convex body,
which is symmetric in respect to $\R^{n}$ ($\R^{n}\subseteq \R^{n+s}$ the
first $n$ coordinates).  For any $a\in\R^{n},b\in \R^{s}$ we define $$\psi(a,b):=
Vol((K-(a,b))^{\circ}).$$
The next two lemmas  deal with two properties of $\psi$. The first property is the inequality $\psi(a,0)\leq\psi(a,b)$, while the  other is
$\psi(-(a,b))=\psi(a,b)$, where $a\in\R^{n},b\in \R^{s}$.

\begin{lemma}\label{Vol}
 Let $K \subseteq \R^{n+s} $ be a centrally symmetric convex body, such that $K$ is symmetric with respect to
 $\R^{n}$. Then, for any $(a,b)\in\R^{n+s} $  $$ Vol((K-(a,0))^{\circ})  \leq Vol((K-(a,b))^{\circ}) .$$
\end{lemma}

\noindent {\bf Proof.}  Define as before $\psi(z):=
Vol((K-z)^{\circ}),\ \ \psi :int(K) \rightarrow \R^{+}$. It is clear
that $\psi$ is strictly convex (Remark \ref{B-infu}), and from the
symmetry of the body $K$ with respect to $\R^{n}$ we have
$Vol((K-(a,b))^{\circ}) = Vol((K-(a,-b))^{\circ})$, where $a\in
\R^{n},b\in\R^{s}$.

Thus, we conclude
\[
\begin{array}{l}
Vol((K-(a,0))^{\circ})= \psi (a,0)
\\
\\
= \psi(\frac{1}{2}(a,b) +
\frac{1}{2}(a,-b))
\\
\\
\leq  \frac{1}{2}\psi((a,b)) + \frac{1}{2}
\psi((a,-b) )
\\
\\
= \psi(a,b)=Vol((K-(a,b))^{\circ}).
\end{array}
\]\qed

\begin{Corollary} \label{corollary2}
Let $K \subseteq \R^{n+s} $ be a centrally symmetric convex body. If
 $(x,y) \in S(K,t)$ then $(x,0) \in S(K,t)$. Indeed, assumption implies
\begin{eqnarray*}
 Vol(K)Vol((K-(x,y))^{\circ}) \leq Vol(B_{2}^{n+s})^{2}t.
\end{eqnarray*}
Hence,  from Lemma \ref{Vol}
\begin{eqnarray*}
 Vol(K)Vol((K-(x,0))^{\circ}) \leq Vol(B_{2}^{n+s})^{2}t,
\end{eqnarray*}
which leads to $(x,0) \in S(K,t)$.
\end{Corollary}

\begin{lemma}\label{newVol} Let $K \subseteq \R^{n+s}$ be a centrally
symmetric convex body. Then  $$Vol((K-z)^{\circ})=Vol((K+z)^{\circ}) \textrm{ for any }z \in
int(K).$$
\end{lemma}

\noindent {\bf Proof.} Since $K=-K$ we have $K+z=-K+z=-(K-z)$,  which implies  $$(K+z)^{\circ}=(-(K-z))^{\circ}=-(K-z)^{\circ}.$$
Hence  $Vol((K+z)^{\circ})=Vol((K-z)^{\circ})$.\qed

 The next lemma deals with the continuity of the {\it Floating body} with
respect to locally uniformly convergence of $s$-concave functions.

\begin{lemma}\label{14}
Let $f: \R^{n} \rightarrow [0,\infty) $ be an even log-concave
function. Define as before $f_{s}(x)= (1+ \frac{logf(x)}{s})_{+}$. Then   $$F(f_{s},\lambda) \rightarrow_{s\rightarrow \infty}
F(f,\lambda) \textrm{ for any  } 0<\lambda< 1/2.$$
\end{lemma}

\noindent {\bf Proof.} Let $\theta\in S^{n-1}$, and let
 $H_{f_{s},\theta}$, $H_{f,\theta}$  be the two hyperplanes
perpendicular to $\theta$ such that
\begin{eqnarray*}
\int_{H_{f_{s},\theta}^{+}}f_{s}(x)dx=\lambda
\int_{\R^{n}}f_{s}(x)dx \ \ \  and \ \ \
\int_{H_{f,\theta}^{+}}f(x)dx=\lambda \int_{\R^{n}}f(x)dx.
\end{eqnarray*}
Since $f_{s}\rightarrow_{s\rightarrow \infty}
f$ locally uniformly we conclude
\begin{eqnarray*}
\int_{H^{+}_{f_{s},\theta}}f_{s}dx=\lambda \int
f_{s}dx\rightarrow_{s\rightarrow \infty} \lambda \int f dx =
\int_{H^{+}_{f,\theta}}f dx.
\end{eqnarray*}
Since $H_{f_{s},\theta}$ and $H_{f,\theta}$ are parallel, we
conclude that $H_{f_{s},\theta}\rightarrow_{s\rightarrow \infty}
H_{f,\theta},$ with respect to the {\it Hausdorff distance} of bodies. \qed

Before we proceed with the proof of Theorem \ref{Thm}, let us prove the following technical lemma.
 \begin{lemma}\label{13}
Let  $ n,s$ be a positive integers. Then
\begin{eqnarray*}
  [ \frac{Vol(B_{2}^{s})^{2}}{Vol(B_{2}^{s+n})^{2}} (\frac{2
\pi}{s})^{n} ]  \rightarrow_{s\rightarrow \infty} 1 \textrm{ from above.}
\end{eqnarray*}
\end{lemma}

\noindent {\bf Proof.}  For $n
\in \N$ even,  $Vol(B_{2}^{n})=\frac{(2\pi)^{\frac{n}{2}}}{n!!}$, while for $n \in \N$ odd  $Vol(B_{2}^{n})=\frac{\pi^{\frac{n-1}{2}}2^{\frac{n+1}{2}}}{n!!}$.  The rest of the proof follows from simple limit calculation.   \qed

 \noindent {\bf Proof of Theorem \ref{Thm}.}  Let $f:\R^{n}\rightarrow [0,\infty)$ be an even log-concave
 function and let $t>0$. As mentioned before,  we first prove Theorem \ref{Thm} with respect to $f_{[t]}$,
  then from the following facts
$$F(f_{[t]},\lambda)\rightarrow_{t\rightarrow_{\infty}}F(f,\lambda) \textrm { and } S(f_{[t]},d)\rightarrow_{t\rightarrow\infty}S(f,d),$$
we will reach the proof.

For $s\in \N$ consider  $f_{[t],s}(x)=(1+
\frac{logf_{[t]}(x)}{s})_{+}^{s}$.
 From Lemma \ref{10} we know
\begin{eqnarray}\label{first}
\{a \in \R^{n}: \int f_{[t],s}  \int
 \L_{s}(f_{[t],s,a})\leq (2\pi)^{n}d \} \rightarrow_{s \rightarrow \infty} S(f_{[t]},
 d).
\end{eqnarray}
\noindent We claim
\begin{eqnarray}\label{second}
  F(f_{[t],s},\lambda)\subseteq \{a \in \R^{n}: \int f_{[t],s}  \int
\L_{s}(f_{[t],s,a})  \leq (2\pi)^{n}d \}.
\end{eqnarray} Indeed, Lemma \ref{13} implies  $[
\frac{Vol(B_{2}^{s})^{2}}{Vol(B_{2}^{s+n})^{2}} (\frac{2
\pi}{s})^{n} ]  \rightarrow_{s\rightarrow\infty} 1$, and $[
\frac{Vol(B_{2}^{s})^{2}}{Vol(B_{2}^{s+n})^{2}} (\frac{2
\pi}{s})^{n} ]\ge1$. Hence
\[
\begin{array}{l}
\big\{a \in \R^{n}: \int f_{[t],s} \int
 \L_{s}(f_{[t],s,a})  \leq (2\pi)^{n}d\big\}
\\
\\
=\bigg\{a \in \R^{n}: Vol(K_{s}(f_{[t],s}))Vol(K_{s}(\L_{s}(f_{[t],s,a})))
\leq [ \frac{Vol(B_{2}^{s})^{2}}{Vol(B_{2}^{s+n})^{2}} (\frac{2
\pi}{s})^{n} ] Vol(B_{2}^{s+n})^{2} d \bigg\}.
\\
\\
\supseteq \bigg\{(a,0) \in \R^{n}\times \R^{s}:
   Vol(K_{s}(f_{[t],s}))Vol(K_{s}(\L_{s}(f_{[t],s,a})))
 \leq  Vol(B_{2}^{s+n})^{2} d
\bigg\}
\\
\\
=\bigg\{(a,0) \in \R^{n}\times \R^{s}:
Vol(K_{s}(f_{[t],s}))Vol(K_{s}(f_{[t],s,a})^{\circ})
 \leq  Vol(B_{2}^{s+n})^{2}d
\bigg\},
\end{array}
\]
\noindent Since $K_{s}(f_{[t],s,a})=K_{s}(f_{[t],s})+ \frac{1}{\sqrt
s}(a,0)$, So we get:
\[
\begin{array}{l}
 \bigg\{(a,0) \in \R^{n}\times \R^{s}:
Vol(K_{s}(f_{[t],s}))Vol(K_{s}(f_{[t],s,a})^{\circ})
 \leq  Vol(B_{2}^{s+n})^{2}d
\bigg\}
\\
\\
= \bigg\{(a,0) \in \R^{n}\times \R^{s}:
Vol(K_{s}(f_{[t],s}))Vol((K_{s}(f_{[t],s})+ \frac{1}{\sqrt s
}(a,0))^{\circ})
 \leq Vol(B_{2}^{s+n})^{2}d
\bigg\}
\\
\\
= \sqrt s \bigg\{(a,0) \in \R^{n}\times \R^{s}:
Vol(K_{s}(f_{[t],s}))Vol((K_{s}(f_{[t],s})+ (a,0))^{\circ})
 \leq Vol(B_{2}^{s+n})^{2}d
\bigg\}.
\end{array}
\]
Since $K_{s}(f_{[t],s})$ is a centrally symmetric convex body, from
Lemma \ref{newVol} we conclude
$$Vol((K_{s}(f_{[t],s})+ (a,0))^{\circ})=Vol((K_{s}(f_{[t],s})-
(a,0))^{\circ}).$$

\noindent From Corollary
\ref{corollary2},   assumption of Theorem \ref{Thm} and Lemma \ref{good}
\[
\begin{array}{l}
\sqrt s \bigg\{(a,0) \in \R^{n}\times \R^{s}:
Vol(K_{s}(f_{[t],s}))Vol((K_{s}(f_{[t],s})- (a,0))^{\circ})
 \leq  Vol(B_{2}^{s+n})^{2}d
\bigg\}
\\
\\
=\sqrt s P_{\R^{n}}[   S(K_{s}(f_{[t],s}), d) ]
\supseteq \sqrt s P_{\R^{n}}[   F (K_{s}(f_{[t],s}),\lambda)]
\\
\\
=F(f_{[t],s},\lambda).
\end{array}
\]
From Lemma
\ref{14} we conclude  $$F(f_{[t],s},\lambda) \rightarrow_{s
\rightarrow \infty} F(f_{[t]},\lambda).$$
From Lemmas \ref{c1},\ref{c2} we complete the proof. \qed

We turn to the proof of Theorem \ref{reverse}.

\noindent {\bf Proof Theorem \ref{reverse}.} We use the same
notation as in Theorem \ref{Thm}. Let $f:\R^{n}\rightarrow
[0,\infty)$ be an even log-concave function and $t>0$.  First we
prove Theorem \ref{reverse} for $f_{[t]}$, then from the facts
 $F(f_{[t]},\lambda)\rightarrow_{t\rightarrow_{\infty}}F(f,\lambda)$ and
  $S(f_{[t]},d)\rightarrow_{t\rightarrow\infty}S(f,d)$ we
will reach the proof.

As before, for $s\in \N$  we consider $f_{[t],s}(x)=(1+
\frac{logf_{[t]}(x)}{s})_{+}^{s}$.  We define
\[
\begin{array}{l}
L_{s}:=\bigg\{ a \in \R^{n}:  \frac{Vol(B_{2}^{s})}{s^{\frac{n}{2}}}\int
f_{[t],s}dx \frac{Vol(B_{2}^{s})}{s^{\frac{n}{2}}}\int
\L_{s}(f_{[t],s,a})dx  \leq   Vol(B_{2}^{n+2})^{2}d\bigg\}
\\
\\
=\bigg\{ a \in \R^{n}:  \int f_{[t],s}dx  \int \L_{s}(f_{[t],s,a})dx \leq
[\frac{Vol(B_{2}^{n+s})^{2}}{Vol(B_{2}^{s})^{2}}
(\frac{s}{2\pi})^{n}]  (2\pi)^{n} d\bigg\}.
\end{array}
\]
In what follows we will show
\begin{eqnarray}\label{Ls convergence}
L_{s}\rightarrow_{s\rightarrow\infty}S(f_{[t]},d).
\end{eqnarray}
\begin{eqnarray}\label{Ls contained}
 L_{s}\subseteq F(f_{[t],s},\lambda).
\end{eqnarray}
Before we continue,   since
$F(f_{[t],s},\lambda)\rightarrow_{s\rightarrow\infty}F(f_{[t]},\lambda)$ (Lemma \ref{14}), then from (\ref{Ls convergence}) and (\ref{Ls contained}) we complete the proof.

We begin with (\ref{Ls convergence}). From lemma \ref{13} we know
$[\frac{Vol(B_{2}^{n+s})^{2}}{Vol(B_{2}^{s})^{2}}
(\frac{s}{2\pi})^{n}] \rightarrow_{s\rightarrow \infty} 1$. Hence
 Lemma \ref{10} implies that $L_{s} \rightarrow_{s\rightarrow \infty} S(f_{[t]},d)$.\\

We continue with (\ref{Ls contained}).   Simple  calculation yields
$$ \frac{Vol(B_{2}^{s})}{s^{\frac{n}{2}}}\int
f_{[t],s}dx=Vol(K_{s}(f_{[t],s}) )$$
 $$
\frac{Vol(B_{2}^{s})}{s^{\frac{n}{2}}}\int \L_{s}(f_{[t],s,a})dx =
Vol(K_{s}(\L_{s}(f_{[t],s,a}))).$$ Thus, we conclude
\begin{eqnarray*}
L_{s}= \{a \in \R^{n}:
Vol(K_{s}(f_{[t],s}))Vol(K_{s}(\L_{s}(f_{[t],s,a}))) \leq
Vol(B_{2}^{n+s})^{2}d\}.
 \end{eqnarray*}
Moreover, from Lemma
\ref{C}
\[
\begin{array}{l}
   \bigg\{(a,0) \in \R^{n}\times \R^{s}: Vol(K_{s}(f_{[t],s}))Vol(K_{s}(\L_{s}(f_{[t],s,a})))
 \leq  Vol(B_{2}^{s+n})^{2} d
\bigg\}
\\
\\
= \bigg\{(a,0) \in \R^{n}\times \R^{s}:
Vol(K_{s}(f_{[t],s}))Vol((K_{s}(f_{[t],s,a}))^{\circ})
 \leq  Vol(B_{2}^{s+n})^{2}d
\bigg\}.
\end{array}
\]

As before, we know  $K_{s}(f_{[t],s,a})=K_{s}(f_{[t],s})+ \frac{1}{\sqrt
s}(a,0)$, which leads to
\[
\begin{array}{l}
 \bigg\{(a,0) \in \R^{n}\times \R^{s}:
Vol(K_{s}(f_{[t],s}))Vol((K_{s}(f_{[t],s,a}))^{\circ})
 \leq  Vol(B_{2}^{s+n})^{2}d
\bigg\}
\\
\\
= \bigg\{(a,0) \in \R^{n}\times \R^{s}:
Vol(K_{s}(f_{[t],s}))Vol((K_{s}(f_{[t],s})+ \frac{1}{\sqrt s
}(a,0))^{\circ})
 \leq Vol(B_{2}^{s+n})^{2}d
\bigg\}
\\
\\
= \sqrt s \bigg\{(a,0) \in \R^{n}\times \R^{s}:
Vol(K_{s}(f_{[t],s}))Vol((K_{s}(f_{[t],s})+ (a,0))^{\circ})
 \leq Vol(B_{2}^{s+n})^{2}d
\bigg\}.
\end{array}
\]

\noindent Since $K_{s}((f_{[t]})_{s})$ is a centrally symmetric convex body,
from Lemma \ref{newVol} we conclude
$$Vol((K_{s}(f_{[t],s})+ (a,0))^{\circ})=Vol((K_{s}(f_{[t],s})-
(a,0))^{\circ}).$$ Thus, from Corollary
\ref{corollary2}, assumption of Theorem \ref{reverse} and Lemma  \ref{good}
\[
\begin{array}{l}
\sqrt s \bigg\{(a,0) \in \R^{n}\times \R^{s}:
Vol(K_{s}(f_{[t],s}))Vol((K_{s}(f_{[t],s})- (a,0))^{\circ})
 \leq  Vol(B_{2}^{s+n})^{2}d
\bigg\}
\\
\\
=\sqrt s P_{\R^{n}}[   S(K_{s}(f_{[t],s}), d) ]
\subseteq
\sqrt s P_{\R^{n}}[   F (K_{s}(f_{[t],s}),\lambda)]
\\
\\
= F(f_{[t],s},\lambda),
\end{array}
\]
which completely shows (\ref{Ls contained}). \qed

\end{document}